\documentclass[lettersize,journal]{IEEEtran}

\usepackage[utf8]{inputenc}
\usepackage{amsmath}
\usepackage{algorithm}
\usepackage{array}
\usepackage[caption=false,font=normalsize,labelfont=sf,textfont=sf]{subfig}
\usepackage{textcomp}
\usepackage{stfloats}
\usepackage{url}
\usepackage{verbatim}
\usepackage{graphicx}
\usepackage{cite}
\usepackage{xspace}

\usepackage{multirow}
\usepackage{xcolor}
\usepackage{colortbl}
\usepackage[mathscr]{euscript}
\usepackage{amssymb,amsbsy}
\usepackage{wrapfig}
\usepackage{fancyhdr}
\usepackage{tikz,pgfplots}
\usepackage{multirow}

\newcommand{\dt}{{\Delta t}}

\newcommand{\bs}[1]{\ensuremath{\boldsymbol{#1}}}

\DeclareMathAlphabet{\mathup}{OT1}{\familydefault}{m}{n}

\definecolor{grassgreen}{RGB}{92,135,39}

\newcommand{\subscript}[2]{{#1}_{{\!\mbox{\rm \tiny #2}}}}
\newcommand{\superscript}[2]{{#1}^{{\!\mbox{\rm \tiny #2}}}}

\renewcommand{\vec}[1]{{\mathchoice
                     {\mbox{\boldmath$\displaystyle{#1}$}}
                     {\mbox{\boldmath$\textstyle{#1}$}}
                     {\mbox{\boldmath$\scriptstyle{#1}$}}
                     {\mbox{\boldmath$\scriptscriptstyle{#1}$}}}}  
\newcommand{\mat}[1]{\bs{#1}}                                      
\newcommand{\diag}{\operatorname{diag}}                            
\newcommand*{\tran}{^{\mkern-1.5mu\mathsf{T}}}                     
\newcommand*{\inv}{^{\mkern-1.5mu\mathsf{-1}}}                     
\newcommand{\inverse}[1]{\Bigl(#1\Bigr)\inv}                       
\newcommand{\Rnum}{\mathbb{R}}                                     
\newcommand{\norm}[1]{\|{#1}\|}                                    
\newcommand{\wnorm}[2]{\norm{#1}_{#2}}                             

\DeclareMathOperator*{\argmin}{arg\,min}  

\newcommand{\ipar}{\vec{\theta}}                   
\newcommand{\iparmap} {\subscript{\ipar}{MAP}}     
\newcommand{\ipartrue}{\subscript{\ipar}{true}}    
\newcommand{\iparprior}{\subscript{\ipar}{prior}}  
\newcommand{\state}{\vec{z}}                       
\newcommand{\simmodel}{\mathcal{S}}                
\newcommand{\tlm}{\mat{S}}                         
\newcommand{\true}{^{{\!\mbox{\rm \tiny true}}}}   

\newcommand{\Nstate}{{\textsc{N}_{\rm state}}}  
\newcommand{\Nparam}{{\textsc{N}_{\rm param}}}  
\newcommand{\Nens}{{\textsc{N}_{\rm ens}}}      
\newcommand{\Nobs}{{\textsc{N}_{\rm obs}}}      
\newcommand{\nobs}{n_{t}}                       
\newcommand{\nobstimes}{\nobs}                  

\newcommand{\noise}{\vec{\zeta}}                   
\newcommand{\covmat}{\mat{\Gamma} }                
\newcommand{\prcov}{\subscript{\covmat}{prior} }   
\newcommand{\postcov}{\subscript{\covmat}{post} }  
\newcommand{\ncov}{\subscript{\covmat}{noise} }    

\newcommand{\GM}[2]{\mathcal{N}\!\left( {#1}, {#2}\right)}  
\newcommand{\probmeasure}{\pi}                              
\newcommand{\like}{\subscript{\probmeasure}{like}}          
\newcommand{\post}{\subscript{\probmeasure}{post}}          
\newcommand{\prior}{\subscript{\probmeasure}{prior}}        

\newcommand{\ff}{{\ensuremath{\mathcal{F}}}}       
\newcommand{\eeta}{\vec{\eta}}

\newcommand{\J}{\mathcal{J}}

\newcommand{\ObsOper}{\mat{O}}                      
\newcommand{\ObsOpertind}[1]{\ObsOper_{\tind{#1}}}  

\newcommand{\tind}[1]{{t}_{#1}}                     
\newcommand{\obs}{\vec{d}}                          

\newcommand{\Exp}[1]{\exp{\left(#1\right)}}         

\usepackage{algpseudocode}

\newcommand{\pyoed}{{PyOED}\xspace }

\newcommand{\gradpower}{{GradPower}\xspace }

\newcommand{\fdvar}{{4DVar}\xspace }

\newcommand{\Adrian}[1]{{\color{olive}Adrian: #1}}

\usepackage[colorinlistoftodos]{todonotes}

\begin{document}

    \title{Centralized calibration of power system dynamic models using variational data assimilation}

    \author{
        Ahmed Attia~\IEEEmembership{Member,~IEEE},
        D. Adrian Maldonado~\IEEEmembership{Member,~IEEE}, 
        Emil Constantinescu~\IEEEmembership{Member,~IEEE},
        Mihai Anitescu~\IEEEmembership{Member,~IEEE}
    \thanks{
        This material is based upon work supported by the
        U.S. Department of Energy, Office of Science, Advanced Scientific
        Computing Research under Contract DE-AC02-06CH11357.
    }
    \thanks{
        Authors are with the Mathematics and Computer Science Division at 
        Argonne National Laboratory, Lemont, Illinois, U.S.A.}
    }%
    
    \markboth{IEEE Transactions on Power Systems}%
    {Optimization-Based Propagation of Uncertainties in Power System Dynamics}
    
    
    \maketitle
    
    \begin{abstract}
    This paper presents a novel centralized, variational data assimilation approach for calibrating transient dynamic models in electrical power systems, focusing on load model parameters. With the increasing importance of inverter-based resources, assessing power systems' dynamic performance under disturbances has become challenging, necessitating robust model calibration methods. The proposed approach expands on previous Bayesian frameworks by establishing a posterior distribution of parameters using an approximation around the maximum a posteriori value. We illustrate the efficacy of our method by generating events of varying intensity, highlighting its ability to capture the systems' evolution accurately and with associated uncertainty estimates. This research improves the precision of dynamic performance assessments in modern power systems, with potential applications in managing uncertainties and optimizing system operations.
    \end{abstract}


    \section{Introduction}
    \label{sec:introduction}

    \IEEEPARstart{E}lectrical power systems are transforming, with inverter-based resources (IBRs) becoming increasingly significant in the energy mix. 
    Because of the rapid dynamics nature of IBRs and their influence on the inertial properties of power systems, assessing the system's dynamic performance under disturbances is becoming progressively more challenging. Given these changes, calibration of transient dynamic models becomes crucial \cite{Zhao2019}. The parameters for these models are often obtained through first principles and experiments, such as for large synchronous generators, or via a blend of engineering judgment and field data recordings.

    One key method for calibrating transient dynamic models is dynamic parameter estimation, also known as \textit{data assimilation} in other fields. This method involves leveraging state estimation techniques such as the Kalman filter or the unscented Kalman filter and augmenting the state vector with unknown or uncertain parameters \cite{Valverde2011, Rouhani2016}. The bulk of the methods in the literature can be classified as \textit{sequential} data assimilation techniques and \textit{decentralized}.

    Data assimilation techniques are typically classified into sequential and variational. Sequential techniques integrate observations as they become available, updating the system based on the current estimate. Conversely, variational techniques consider an entire time window of observations to minimize the discrepancy between prediction and observation. While sequential techniques tend to be more straightforward to implement and computationally efficient, variational methods offer more robustness, particularly when model error and nonlinearity are present \cite{Constantinescu_A2007b,attia2015hmcfilter,attia2015hmcsampling,attia2015hmcsmoother}. While empirically filtering methods can be nonetheless competitive even in nonlinear settings for the proper parameter choices, we will restrict our investigation to variational approaches (which have maximum likelihood flavor) for the reasons expressed in this paragraph.  
    
    A further feature of the calibration problem involves uncertainty quantification. Some calibration methods aim to identify model parameters that align with specific data. However, determining the uncertainty associated with such estimates is often beneficial. For this, Bayesian approaches are useful, as they provide a statistical framework for incorporating prior knowledge and observed data to produce a probabilistic description of the unknown parameters.

    A secondary classification is centralized and decentralized in a sense employed by\cite{Zhao2019}. Decentralized approaches examine small subsystems, seeking local optimal solutions to the calibration problem. For example, methods for calibrating synchronous generator models using phasor measurement units (PMUs)  at the terminal are commonly found in the literature. Since the voltage and current phasors fully determine the synchronous generator model behavior, measurements from other buses may not add any necessary information to the estimation. However, when calibrating aggregated load models or IBR models, which often inaccurately represent actual behavior, PMU measurements may miss higher frequency dynamics impacting several inverters in an area. In such cases, a decentralized calibration approach might produce optimal parameters for a specific event but fail to generalize to others. Moreover, when considering the uncertainty representation, decentralized approaches cannot describe the correlation between the parameter uncertainties, thus missing important response features. Such situations can be mitigated by adopting a centralized approach, considering the entire system model. We emphasize that here we used the words centralized and decentralized relative to the model used in the data assimilation; centralized (global) models can still be solved in principle with parallel, decentralized algorithms. 

        Our work introduces a centralized, variational data assimilation approach for model calibration. We expand the Bayesian framework work in \cite{petra2016bayesian} and the recent developments in \cite{nagi2021bayesian,saric2018data}. Petra et al. 2016 \cite{petra2016bayesian} demonstrate a Bayesian approach for parameter inversion that uses a Laplace approximation around the maximum a posteriori (MAP) point to characterize the posterior, which does not rely on Gaussian priors and likelihoods. Nagi et al. 2022 \cite{nagi2021bayesian} propose a Bayesian parameter inversion that leverages trajectory sensitivities for the linearized power grid model. The linear aspect allows one to find the exact analytical solution resulting from computed Gaussian posteriors using conjugate distributions for the prior and likelihood. The authors use this solution  to achieve high-dimensional scalability. Alternative inference approaches based on polynomial chaos expansion (PCE) have also been used to compute the posterior distribution in Bayesian parameter inversion \cite{xu2019response,xu2019bayesian,petra2016bayesian}. While PCE is robust to the prior and likelihood distributions, high parameter dimension inferences are problematic and require special attention. We note that all these works, including the one introduced here, assume a centralized parameter estimation; however, conditional independence of load characteristics at every bus may allow in principle independent calibration at every bus if the bus is instrumented to measure the power flow. We note, however, that, as discussed above, this may result in wider confidence intervals compared with centralized approaches. 
        
        In this study we focus on the estimation of parameters of load models. We assume that the load model composition is fixed, once calibrated during the prediction window. This aspect can be extended by ($i$) training continuously or at fixed intervals or ($ii$) using a data-driven approach in which the load parameters are aggregated across multiple inference windows. We also assume that the initial disturbance is fully characterized. One potential approach to relax this assumption is to include these uncertainties in the inference model. Nevertheless, this study provides accurate predictions in this context, and these results can be used to extend these ideas to more complex situations.  

This paper's main contribution is developing a posterior distribution of the parameters by approximating the probability density function around the maximum posterior value of the parameter, along with strategies to accommodate compact supports. This will allow us to gauge the benefit of linearizing the system at the mode of the Bayesian posterior, as opposed to linearizing the system at a reference point as carried out in  \cite{nagi2021bayesian}, while in our case facing the additional difficulty of carrying out the nonlinear maximization. 
        
In Section \ref{sec:models} we describe the dynamic model, the load model, and the numerical integration method. In Section \ref{sec:Bayesian_inversion} we present the mathematical formulation of the calibration algorithm, which is formulated as a Bayesian inversion problem. In Section \ref{sec:numerical_results} we test our algorithm by simulating diverse fault scenarios and using the resulting information as input data to calibrate load parameters. In addition, we show how our algorithm outperforms the state of the art.

    \section{Dynamics and measurement models}
    \label{sec:models}
        Our problem consists of calibrating the parameters of a set of load models 
        given a collection of recorded events. 
        We assume a differential-algebraic equations (DAE) model of 
        the positive-sequence dynamics. 
        Recorded events are obtained by using PMUs 
        with exhaustive or sparse measurements.
    
        \subsection{Power system dynamic model}\label{sec:PoweSystemModel}
            The dynamics of the power system are modeled by using the following
            DAE system:
                \begin{subequations}\label{eqn:powersystem}
                  \begin{align}
                    &\dot{\vec x} = f(t,\vec x,\vec y;\ipar)\,,\label{eqn:ode}\\
                    &\vec 0 = g(t,\vec x,\vec y;\ipar)\,,\label{eqn:algeqs} \\
                    &\vec x(0) = \vec x_0\,,  \quad  \vec y(0) = \vec y_0\,. \label{eqn:ic}
                  \end{align}
                \end{subequations}

                Here $\vec x \in \mathbf{R}^n$ represents dynamic variables such as rotor
                angle or generator speed; $\vec y \in \mathbf{R}^m$ represents the algebraic variables
                such as the bus voltages and line currents; $\vec x_0$ and $\vec y_0$ are
                the initial conditions; $t$ represents time; 
                and $\ipar \in \mathbf{R}^{\Nparam}$ denotes the model parameters.
                The right-hand side $f(\cdot)$ in~\eqref{eqn:ode} is in general a nonlinear function 
                that models the dynamics of the system, and $g(\cdot)$ in~\eqref{eqn:algeqs} is 
                a set of algebraic equations modeling the passive network of the power system.
                In general, $f(\cdot)$ models dynamic devices such as generators, governors, exciters, 
                and dynamical loads, whereas $g(\cdot)$ models the network balance equations and 
                the current injections of each device connected to it. 
                
                In this work we focus on the calibration of passive load models of the form
                \begin{align}
                \begin{split}
                \label{eq:load}
                p &= \theta_i \left(  \frac{v}{v_0} \right)^2 p_0 + (1 - 
                \theta_i) p_0  \,, \\
                q &= \theta_i \left(  \frac{v}{v_0} \right)^2 q_0 + (1 - 
                \theta_i) q_0  \,,            
                \end{split}
                \end{align}
                where, for a given bus, $(p, q)$ are the instant active and reactive power consumption, respectively;
                $(p_0, q_0)$ are the base active and reactive power consumption, respectively; 
                $v$ is the voltage;
                $v_0$ is the steady-state voltage magnitude;
                and $\theta_i$ is a component of $\ipar$ that determines the mixture between constant power and constant impedance. 
                While quantities such as $p_0$ and $v_0$ are parameters, for clarity, we  include in $\ipar$ only those parameters that will be calibrated (inferred).
                To further abstract the DAE model, we write the system~\eqref{eqn:powersystem} in the following form:
                \begin{subequations}\label{eqn:powersystem_agrregated}
                    \begin{align}
                        &\mat{M} \dot{\state} = h(t, \state; \ipar) \,, \\
                        & \state(0) = \state_{0} \,,
                    \end{align}
                \end{subequations}
                where $\state \in \Rnum^{\Nstate} $ denotes the model state that 
                aggregates $\vec{x}, \vec{y}$, that is, $\Nstate = m + n$, 
                and $\mat{M} \in \Rnum^{\Nstate \times \Nstate}$
                is the \textit{mass} matrix
                \begin{equation*}
                    \mat{M} = \begin{bmatrix}
                    \mat{I} & \mat{0} \\
                    \mat{0} & \mat{0}
                    \end{bmatrix} \,, \quad
                    \mat{I} \in \Rnum^{n \times n} \,.
                \end{equation*}

                We chose to calibrate a load model because this is one of the less understood 
                and most difficult to measure parts of the power system. 
                Furthermore, while these models are simple, 
                the paucity of \textit{informative} data makes overfitting 
                a real danger \cite{Maldonado2017}.
                Our framework, however, can handle any DAE model calibration
                and is trivial to adapt to calibrate parameters of generators, 
                inverter-based resources, or any other element with unknown parameters.

        \subsection{Forward integration and trajectory sensitivities}
            \label{subsec:forwrad_model}
            The DAE model introduced in equations \eqref{eqn:powersystem} 
            a--c can be integrated by different numerical methods. One such method is backward Euler, which we will use in our work. Given \eqref{eqn:powersystem_agrregated}, we solve the system of nonlinear equations
            \begin{equation}\label{eqn:powersystem_beuler}
                \mat{M} z_{k + 1} - \mat{N}h(t, z_{k + 1}; \ipar) = \mat{M} z_k  \,,
            \end{equation}
            where $z_{k+1}$ is the unknown state at $k+1$ step of the state trajectory, $N$ is defined as
            \begin{equation*}
                    \mat{N} = \begin{bmatrix}
                    \dt \mat{I}_n & 0 \\
                    0 & \mat{I}_m
                    \end{bmatrix} \,, \quad
                    \mat{I}_n \in \mathbf{R}^{n \times n}, \mat{I}_m \in \mathbf{R}^{m \times m} \,,
            \end{equation*}
            and $\dt$ is the discretization time step size. 
            At the same time, we solve another system of linear, time-variant ODEs that describe the evolution of the trajectory sensitivities.
            This is the tangent linear model (TLM), and we direct the reader to our previous work that includes a description of the model~\cite{maldonado2022trust}.

        \subsection{Measurement model}
            We assume measurements are obtained by PMUs situated at diverse bus locations and obtaining voltage, current injection, and current measurements. Following~\cite{Zhou2015}, we model the PMU measurements as the phasor voltage and current quantities of the DAE plus an additive noise
            that represents measurement uncertainty.
            
            In this work we measure the voltage of all buses in the IEEE 39-bus power grid; 
            however, as is the case for all other state space models~\cite{crassidis2011optimal},
            our framework can be used to experiment with various measurement scenarios
            (e.g., measurements at a subset of buses) at various time intervals
            and measurements of different quantities.

    \section{Parameter Inference and UQ}
    \label{sec:Bayesian_inversion}
        In this study we consider the inference parameter $\ipar$ to be 
        the load parameters $\ipar\in\Rnum^{\Nparam}$.
        We also assume that the initial state $\state_0$ has 
        already been inferred from data and is known. 
        Our focus is on characterizing the state of the network and providing 
        accurate projections of the system dynamics in case of contingencies.
        Nevertheless, initial conditions can also be considered unknown; 
        and the framework introduced herein naturally extends to such cases
        where both the calibration parameter $\ipar$ and initial condition $\state_0$
        can be simultaneously inferred from noisy data.
        In this section we describe the approach proposed for parameter $\ipar$ 
        identification and UQ.

        \subsection{Bayesian inversion}
            Assume that we have measurements (forward problem) of 
            a dynamic system that can be modeled by 
            an additive Gaussian noise model,
            \begin{equation}\label{eqn:noise-model}
                \obs = \ff(\ipar) + \vec{\eta}, \quad \eeta \sim
                \GM{\vec{0}}{\ncov}\,,
            \end{equation}
            where $\ncov \in \Rnum^{\Nobs \times \Nobs}$ is the measurement noise covariance
            matrix and $\ff(\cdot)$ is a (generally nonlinear) operator mapping model
            parameters $\ipar$ to observations $\obs$. 
            In this study the observations correspond to the voltage and current phasors measured by PMUs
            at predefined observation time instances.

            We follow a Bayesian formulation that poses the parameter identification problem 
            as a problem of statistical inference over the parameter space. 
            The solution of the resulting Bayesian inverse problem is 
            a {\it posterior probability density function} (PDF) $\post(\ipar|\obs)$ 
            formulated by applying Bayes' rule,
            \begin{align} \label{eqn:Bayes}
                \post(\ipar|\obs)
                \! \propto \like(\obs | \ipar ) \; \prior(\ipar) \,,
            \end{align}
            where $\propto$ implies dropping the PDF normalization constant,
            which generally requires extensive computations and may mandate 
            Monte Carlo approximations that are extremely challenging especially 
            in high dimensions and large-scale applications.
            
            The likelihood model $\like(\obs | \ipar )$ is derived from the forward problem~\eqref{eqn:noise-model}
            and is thus assumed to be Gaussian $\like(\obs | \ipar)\sim \GM{\ff(\ipar)}{\ncov}$.
            Since the model $\ff$ is nonlinear, the posterior~\eqref{eqn:Bayes} is non-Gaussian 
            and is generally intractable.
            The prior $ \prior(\ipar)$ encapsulates the knowledge about 
            the inference parameter $\ipar$ before data acquisition and assimilation
            and can also encode physics or regularization constraints.
            
            We approximate the posterior distribution by following a Laplacian approach in which 
            the posterior is approximated by a Gaussian that captures both the central tendency of the posterior 
            and covariances quantifying posterior uncertainties. 
            Specifically, the posterior central tendency is approached by seeking the 
            MAP estimate, and the posterior covariance matrix is then 
            approximated around that MAP estimate.
            The proposed approach is discussed in Section~\ref{subsec:FDVar}.

        \subsection{Inverse problem and UQ}
        \label{subsec:FDVar}
            Here we describe the elements of the inverse problems. 
            We describe the likelihood, the prior, and the posterior obtained by~\eqref{eqn:Bayes}, 
            followed by the proposed approach.  
            
            \paragraph{The likelihood}
            \textit{Parameter-to-observable} map $\ff(\ipar)$ evaluation
            requires solving the DAE system that models the dynamics of 
            the power grid, followed by extraction of the PMU measurements from the DAE solution at 
            observation time instances $\tind{1}, \ldots, \tind{\nobstimes}$.
            Let the solution operator $\simmodel$ denote forward 
            integration (simulation) of the state that is 
            the evaluation of the DAE for a given
            model parameter $\ipar$, and the fixed initial model state $\state_0$.
            Specifically, the model state $\state_k$ at time instance $\tind{k}$ is defined as
            \begin{equation}\label{eqn:discrete_forward_model}
                \state_{k} 
                    := \simmodel_{\tind{0} \rightarrow {\tind{k}}}(\ipar, \state_{0})
                    = \simmodel_{0, k}(\ipar, \state_{0}) 
                    = \simmodel_{0, k}(\ipar )
                    \,,
            \end{equation}
            where we dropped the initial state $\state_{0}$ since it is assumed 
            to be known and kept fixed.
            Following~\eqref{eqn:noise-model}, the observations collected at time 
            instance $\tind{k}$ are related to the model state $\state_{k}$,
            by the additive Gaussian noise model,
            \begin{equation}\label{eqn:discrete-noise-model}
                \obs_{k} = \ObsOper_{\tind{k}} (\state_{k}) + \noise\,, 
                    \quad \noise \sim \GM{\vec{0}}{\mat{R}}\,,
            \end{equation}
            where $\ObsOpertind{k}$ is an observation operator that maps 
            the model state $\state_{k}$ onto the observation space at 
            observation time instance $\tind{k}$. 
            Here $\noise$ models observation  
            errors at a given time instance, which is 
            generally assumed to be unbiased with covariance matrix $\mat{R}$.
        In this work we will assume that the observation operator is time-independent, and thus
            $\ObsOper_{\tind{k}}=\ObsOper,\, \forall\, k$.
            The observational errors in this problem are temporally uncorrelated,
            and thus the likelihood is given by
            \begin{equation}\label{eqn:Gaussian_Likelihood}
                 \like(\obs | \ipar ) \propto 
                    \Exp{ 
                        -\frac{1}{2} 
                            \sum_{k=1}^{\nobstimes} 
                                { \wnorm{ \obs_{k}- \ObsOper\left(\simmodel_{0, k}(\ipar)\right)}{\mat{R}\inv}^2 }
                        }  \,,
            \end{equation}
            where the weighted norm is given by $ \wnorm{\vec{x}}{\mat{A}}^2 =\vec{x}\tran \mat{A}\vec{x} $ 
            for a vector $\vec{x}$ and a matrix $\mat{A}$ of conformable sizes.

            \paragraph{The prior}
                The load parameters are constrained to a subset 
                $\ipar\in \Omega \subset \Rnum^\Nparam$, which accounts for the compact support
                of the parameters (typically $\ipar \in [0,1]^\Nparam$).
                Thus, we employ a truncated multivariate Gaussian prior 
                with mean $\iparprior$ and covariance $\prcov$:
                \begin{equation}\label{eqn:truncated_Gaussian_prior}
                    \begin{aligned}
                        \prior(\ipar) 
                            &= 
                            \frac{ \exp{ \left(-\frac{1}{2}  \wnorm{\ipar-\iparprior}{\prcov\inv}^2 \right)} }
                                 { \int_{\Omega}{ \exp{ \left(-\frac{1}{2} \wnorm{\ipar-\iparprior}{\prcov\inv}^2  \right)}} d\ipar }  \\
                            & \propto \exp{\left(-\frac{1}{2} \wnorm{\ipar-\iparprior}{\prcov\inv}^2  \right)} \mathtt{1}_{\Omega}(\ipar) \,, \\
                        \mathtt{1}_{\Omega} (\ipar)
                            &:= \begin{cases}
                                1 \,; & \ipar \in \Omega \\
                                0 \,; & \ipar \notin \Omega \,.
                              \end{cases}
                    \end{aligned}
                \end{equation}
                %
                %

            \paragraph{The posterior}
            Applying Bayes' rule~\eqref{eqn:Bayes} to the Gaussian likelihood~\eqref{eqn:Gaussian_Likelihood} 
            and the truncated Gaussian prior~\eqref{eqn:truncated_Gaussian_prior} yields a truncated posterior 
            \begin{align}\label{eqn:posterior}
                \post(\ipar|\obs) &\propto \exp( - \J(\ipar)) \, \mathtt{1}_{\Omega}(\ipar) \,, \\
                \J(\ipar) 
                    &= \frac{1}{2} \left(
                    \sum_{k=1}^{\nobstimes} 
                        { \wnorm{ \obs_{k} \!-\! \ObsOper\left(\simmodel_{0, k}(\ipar)\right)}{\mat{R}\inv}^2 }
                    + \wnorm{ \ipar \!-\! \iparprior }{\prcov\inv}^2 
                    \right) \,.
            \end{align}

            Because of the nonlinearity of the simulation model $\simmodel$, 
            the posterior~\eqref{eqn:posterior} is non-Gaussian.
            It can, however, be approximated  by a truncated Gaussian 
            $\GM{\iparmap}{\postcov}\, \mathtt{1}_{\Omega}(\ipar)$, 
            where $\iparmap$ is the MAP estimate and $\postcov$
            is an approximate posterior covariance matrix. 
            The MAP estimate $\iparmap$ is
            obtained by maximizing the negative log of the posterior, 
            that is, by solving the following four-dimensional variational (\fdvar) 
            optimization problem (see, e.g.,~\cite{attia2015hmcsmoother,petra2016bayesian}):
            \begin{equation}\label{eqn:MAP_Optimization}
                \iparmap := \argmin_{\ipar \in \Omega} \, \J(\ipar) \,.
            \end{equation}

            Gradient-based constrained numerical optimization algorithms, such as the limited-memory Broyden--Fletcher--Goldfarb-Shanno (L-BFGS-B) algorithm, 
            can be used to numerically solve~\eqref{eqn:MAP_Optimization}.
            The gradient of the objective $\J$ in~\eqref{eqn:MAP_Optimization} is given by
            \begin{equation}\label{eqn:FDVar_gradient}
                \nabla_{\ipar} \J(\ipar) = 
                    \sum_{k=1}^{\nobstimes}
                        {\tlm_{0, k}\tran \ObsOper\tran \mat{R}\inv \left(\obs_{k} \!-\! \simmodel_{0, k}(\ipar) \right)}
                    + \prcov\inv \left( \ipar\!-\!\iparprior \right) \,,
            \end{equation}
            where $\tlm_{0, k}$ is the tangent linear (forward sensitivities), 
            that is the Jacobian of the simulation model 
            $\tlm_{0, k} = \frac{\partial \simmodel_{0, k}}{\partial \ipar}$. 
            The gradient is projected onto the feasible parameter domain $\Omega$ 
            and is then used in a gradient-descent approach to seek a local 
            optimum of~\eqref{eqn:MAP_Optimization}.
            Here we utilized the fact that the observation operator is linear (by definition),
            and thus observation sensitivities are given by $\partial{\ObsOper}=\ObsOper$.
            
            The posterior covariance $\postcov$ is approximated by the linearization of the 
            nonlinear solution model and is given by
            \begin{equation} \label{eqn:Truncated_Gaussian_Posterior_Covariance}
                \postcov = \inverse{ 
                    \sum_{k=1}^{\nobstimes}
                        {\tlm_{0, k}\tran \ObsOper\tran \mat{R}\inv \ObsOper\tlm_{0, k}}
                        + \prcov 
                    }
                    \,,
            \end{equation}
            where the TLM $\tlm_{0, k}$ is evaluated at the MAP estimate $\iparmap$.
            Uncertainty cones (confidence intervals) of the expected observations
            can be obtained by applying the forward operator $\ff$ to samples collected 
            from the truncated Gaussian posterior described above.
            Because sampling a Gaussian requires shifting by the central tendency (MAP) 
            and scaling by the lower Cholesky factor of the covariance matrix, one needs 
            to factorize the posterior covariance~\eqref{eqn:Truncated_Gaussian_Posterior_Covariance}.
            If the parameter dimension is small, one can apply Cholesky factorizationto~\eqref{eqn:Truncated_Gaussian_Posterior_Covariance} directly.
            If the parameter dimension is relatively high, however, one can apply the factorization 
            to the inverse of the posterior covariance matrix, that is, 
            $\mat{L}\mat{L}\tran = \sum_{k=1}^{\nobstimes}
                {\tlm_{0, k}\tran \ObsOper\tran \mat{R}\inv \ObsOper\tlm_{0, k}}
                + \prcov $, and then obtain $\postcov=\superscript{\mat{L}}{-T}\mat{L}\inv
            $,
            where $\mat{L}$ is a lower (Cholesky factor) triangular matrix.

            An algorithmic description of the steps required for evaluating the MAP estimate 
            $\iparmap$ and formulating both state and observations uncertainty cones is
            given by Algorithm~\ref{alg:4DVar_Laplace_Approximation}.
            \begin{algorithm}[htbp!]
                \caption{
                    \fdvar with truncated Gaussian posterior 
                }
                \label{alg:4DVar_Laplace_Approximation}
                \begin{algorithmic}[1]
                
                    \Require{$\iparprior$, $\prcov$, 
                        $\{\obs_{1},\ldots,\obs_{\nobstimes}\}$, $\mat{R}$,
                        $\ObsOper$, 
                        $\simmodel$,
                        $\{\tind{1}, \ldots, \tind{\nobstimes}\}$,
                        $\Nens$
                    }
                    \Ensure{$\iparmap$, $\subscript{\obs}{pred}$, $\subscript{\obs}{pred}\mp 2\sigma$}

                    \State $\iparmap \gets $ solve~\eqref{eqn:MAP_Optimization} 
                        \Comment{gradient-based using~\eqref{eqn:FDVar_gradient}}

                    \State $\tlm_{0, k} \gets $  TLM at $\iparmap$ for $k=1,\ldots,\nobstimes$
                        \Comment{as in~\eqref{eqn:powersystem_beuler}, \cite{maldonado2022trust}}
                    
                    \State Compute the Cholesky factor $\mat{L}$:
                    
                        $
                            \mat{L}\mat{L}\tran= \sum_{k=1}^{\nobstimes}
                            {\tlm_{0, k}\tran \ObsOper\tran \mat{R}\inv \ObsOper\tlm_{0, k}}
                            + \prcov
                        $

                    \State $\vec{\sigma} \gets \sqrt{\diag{\postcov}}$ 
                        \Comment{$\postcov$ given by~\eqref{eqn:Truncated_Gaussian_Posterior_Covariance}, $\sqrt{\cdot}$ is elementwise}

                    \State $i \gets 0$
                    \While{$i<\Nens$}
                        \State Sample a standard normal vector $\vec{z}$ of size $\Nparam$
                        \State $\ipar \gets \iparmap + \mat{L} \vec{z}$
                        \If{$\ipar \in \Omega $ and $\iparmap-2\vec{\sigma}\leq\ipar\leq\iparmap=2\vec{\sigma}$}
                            \State Predicted states:
                                
                                $\qquad
                                \state_{k}^{(i)} \gets \simmodel_{0, k}(\ipar, \state_{0}),\, k=1, \ldots, \nobstimes$
                                    \Comment{Use~\eqref{eqn:discrete_forward_model}}
                                 
                            \State Predicted observations:
                                
                                $\qquad
                                    \obs_{k}^{(i)} \gets \ObsOper(\state_{k}^{(i)}),\, k=1, \ldots, \nobstimes
                                $
                            \State $ i \gets i+1 $
                        \EndIf
                    \EndWhile

                    \State State UQ (standard deviations): $\forall k=1, \ldots, \nobstimes$,
                        
                        $
                            \vec{\sigma}_{k} 
                                = \sqrt{
                                    \frac{
                                        \sum_{i=1}^{\Nens} \left( \state_{k}^{(i)} - \overline{\state_{k}}  \right)^2
                                      }{\Nens-1}
                                  }\,,\quad 
                            \overline{\state_{k}} = \frac{1}{\Nens} \sum_{i=1}^{\Nens} \state_{k}^{(i)}
                        $

                    \State Observation UQ (standard deviations): $\forall i=1, \ldots, \nobstimes$,
                        
                        $
                            \vec{\sigma}_{k} 
                                = \sqrt{
                                    \frac{
                                        \sum_{i=1}^{\Nens} \left( \obs_{k}^{(i)} - \overline{\obs_{k}}  \right)^2
                                      }{\Nens-1}
                                  }\,,\quad 
                            \overline{\obs_{k}} = \frac{1}{\Nens} \sum_{i=1}^{\Nens} \obs_{k}^{(i)}
                        $
                    
                \end{algorithmic}
            \end{algorithm}

    \section{Numerical Results}
    \label{sec:numerical_results}
    In this section we perform a series of experiments on the IEEE 39 bus system. We simulate the IEEE 39 bus system with a positive-sequence model, round rotor synchronous generators (Sauer--Pai model), and turbine governors (IEESGO model). 
    The load, as in \eqref{eq:load}, is modeled as a mixture of constant current and constant power load.
    Event data is obtained by applying three-phase-to-ground faults to selected buses of varying impedance.
    
    \subsection{Experimental setup} 
    \label{subsec:setup}
        All numerical experiments are carried out by using 
        \pyoed~\cite{attia2023pyoed,attia2023pyoedRepo} and \gradpower~\cite{maldonado2023uqgridRepo}.
        Specifically, \gradpower provides forward and adjoint (transposed TLM) evaluations of the simulation model~\eqref{eqn:powersystem}, 
        while \pyoed orchestrates the whole DA process, 
        including evaluation of the parameter-to-observable map $\ff$;
        observations, observation operator, and observation error model; 
        prior and posterior evaluation, quantifying state and observation uncertainties; 
        and visualization.

        \paragraph{Simulation model, TLM, and transients}
            The simulation model and its associated TLM (forward sensitivities) are derived from 
            the standard IEEE 39-bus model and discretization outlined in \S\ref{sec:PoweSystemModel},
            incorporating parameterized loads. By applying faults of varying impedance levels at different buses
            (one combination for each scenario) and clearing them after two cycles, a set of system trajectories is obtained. 
            The system response, obtained using these trajectories, reveals the load composition at  
            observed and unobserved buses. Lower impedance levels cause larger disturbances, yielding more insights into 
            the load composition. Indeed, low-impedance faults induce a large transient that better \textit{explores} 
            the (nonlinear) dynamics of the system. This data is then integrated into the model, producing a posterior 
            distribution of the parameters. The process is detailed below.

        \paragraph{Ground truth and prior}
            We fix a set of values of the parameters as a ground truth, which is used for performance evaluation and for generating synthetic 
            observations for verification. The ground truth $\ipartrue$, that is, the true value of the parameter $\ipar$,
            is set to $\ipartrue=(0.9, 0.9, \ldots, 0.9)\tran\in\Rnum^{\Nparam}$.
            In our example, the number of parameters is set to $\Nparam=19$.
            
            The truncated Gaussian model~\eqref{eqn:truncated_Gaussian_prior} is used to define the prior, 
            where we set the prior mean 
            $\iparprior = (0.5,\ldots, 0.5)\tran \in \Rnum^{\Nparam}$ 
            and the prior covariance $\prcov = 0.01 \mat{I}$,
            where $\mat{I} \in \Rnum^{\Nparam \times \Nparam}$ is the identity matrix. 
            The model parameter $\ipar$ is constrained to $\Omega=[0, 1]^{\Nparam}$. 
            %

        \paragraph{Observations and observation operator}
            The model is simulated over a time interval $[\tind{0}, \tind{f}]=[0, 1]$ seconds, with observations collected 
            at a temporal frequency of $30$ observations per second. 
            Specifically, observations are collected (and assimilated) at time instances
            $\tind{0}+ k \Delta t,$ with $\tind{0}=0.0084,\,k=1, 2, \ldots,$ and $\Delta t=1/29.75$, 
            resulting in $30$ observation time instances $\tind{1}=0.0084, \tind{2}=0.0420, \ldots,\tind{f}=\tind{\nobstimes}=\tind{30}=0.9832$.
            %
            
            The true parameter $\ipartrue$ is propagated forward to generate a reference state 
            and observation trajectories.
            The ground truth of the observations is obtained by applying the forward operator
            $\ff$ to the true parameter $\ipartrue$ at the assimilation time instances.
            %

            The observation operator $\ObsOper$ maps the state vector into the voltage and current 
            phasors as measured by PMUs. Since both voltage and current are part of the  system state, 
            the observation operator here is a linear transformation, and thus the tangent linear (i.e., the Jacobian) 
            of the observation operator $\ObsOper$ satisfies $\partial \ObsOper = \ObsOper$.
            %

        \paragraph{Observation noise and synthetic observations}
            The observation error covariance matrix $\mat{R}$ in~\eqref{eqn:discrete-noise-model} models
            the measurement error covariance at a given time instance $\tind{k}$. 
            Here, $\mat{R}$ is taken as a diagonal matrix similar to previous studies~\cite{nagi2021bayesian},
            which implies that the measurement errors are independent.
            In our framework, however, we can easily incorporate correlations either between voltage-current or between PMUs in different buses. 
            Synthetic observations are created by adding random perturbations 
            (sampled from the observation noise model $\GM{\vec{0}}{\mat{R}}$) 
            to the observed ground truth at the predefined observation/assimilation
            time instances $\tind{0},\ldots,\tind{\nobstimes}$.

        \paragraph{Data assimilation} 
            The initial guess of the optimization procedure for finding the MAP estimate is set to the prior mean $\iparprior$.
            The MAP estimate $\iparmap$ is found by solving~\eqref{eqn:MAP_Optimization} using the selected optimization routine discussed below.
            Once the MAP estimate $\iparmap$ is found, a Gaussian (Laplacian) approximation of the posterior is formulated by using~\eqref{eqn:Truncated_Gaussian_Posterior_Covariance}.
            Note that we do not actually need to construct the posterior covariance matrix $\postcov$, but we need its Cholesky factor for
            posterior sampling as described by Algorithm~\ref{alg:4DVar_Laplace_Approximation}.

            The variational optimization problem~\eqref{eqn:MAP_Optimization} is solved numerically 
            by using an L-BFGS routine~\cite{zhu1997algorithm} in order to enforce the bound constraints $\ipar\in\Omega$; 
            that is, $\ipar_i \in [0, 1]\, \forall\, i=1,\ldots,\Nparam.$ 
            The optimization step size (learning rate) is optimized by using  a standard line-search approach. 
            Here we use a \pyoed optimization routine that employs an L-BFGS-B implementation
            provided by \textsc{SciPy}~\cite{2020SciPy-NMeth}.
            The maximum number of iterations is set to $30$, and the algorithm 
            terminates when the magnitude of the projected gradient is small 
            enough, that is, when a local optimum is found. 
            This is achieved by setting the \textsc{pgtol} parameter to $1e-5$.

    \subsection{Comprehensive numerical study}
        With the experimental setup described above, we perform three sets of experiments.
        First, we analyze the performance of the proposed approach for parameter identification 
        for a given contingency, that is, for a fixed fault impedance level and faulty bus; see~\ref{subsubsec:inversion_results}.
        Second, in~\ref{subsubsec:predectivity} we test the accuracy of the prediction against many contingencies
        given a posterior obtained for a specific contingency.
        Third, in~\ref{subsubsec:tlm_results} we test the proposed approach against the commonly employed linearization approach to UQ.

    \subsubsection{Bayesian inversion results}
    \label{subsubsec:inversion_results}
        In this section we analyze the inversion accuracy of the approach described by Algorithm~\ref{alg:4DVar_Laplace_Approximation}.
        Specifically, we solve the inverse problem by applying Algorithm~\ref{alg:4DVar_Laplace_Approximation} for multiple contingencies. 
        First, we fix the impedance level and change the location where the fault happens. Second, we fix the bus at which the fault occurs and change the impedance level.
        To assess the calibration performance, we plot the root-mean-squared error (RMSE) results of the observations predicted by using the MAP estimate of the parameters, 
        along with posterior uncertainty cones constructed by sampling the truncated posterior and evaluating observations of 
        the posterior trajectory (model state trajectory obtained by propagating the MAP point over the assimilation window).
        The RMSE evaluated at the $k$th time instance $\tind{k}$ for a state/observation vector $\vec{x}_{k}\in\Rnum{n}$, 
        compared with ground truth $\vec{x}_{k}\true$, is given by
        \begin{equation} \label{eqn:RMSE}
            RMSE = \sqrt{
                \frac{1}{n} \left(
                    \vec{x} - \vec{x}\true
                \right)^2
            }  \,.
        \end{equation}
        %
        
\begin{figure}[htbp!]
            \centering
            \includegraphics[trim={0 50pt 0 0},clip,width=0.75\linewidth]{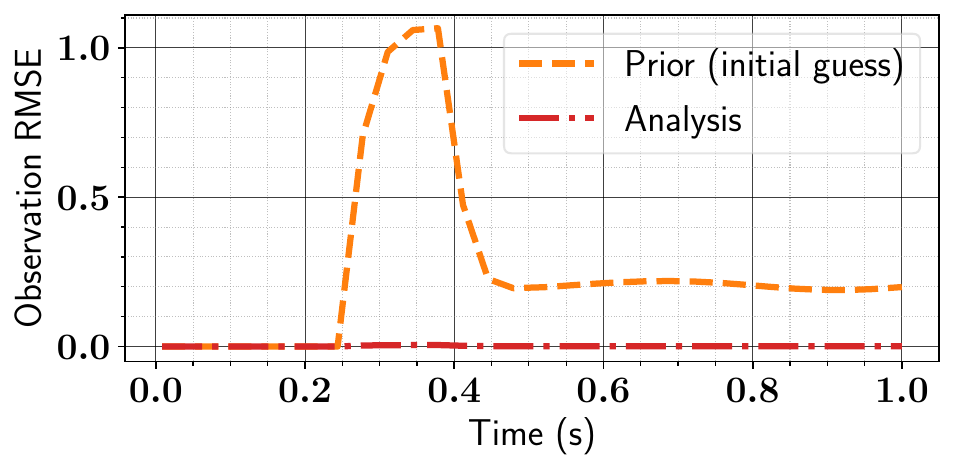}
            \includegraphics[trim={0 50pt 0 0},clip,width=0.75\linewidth]{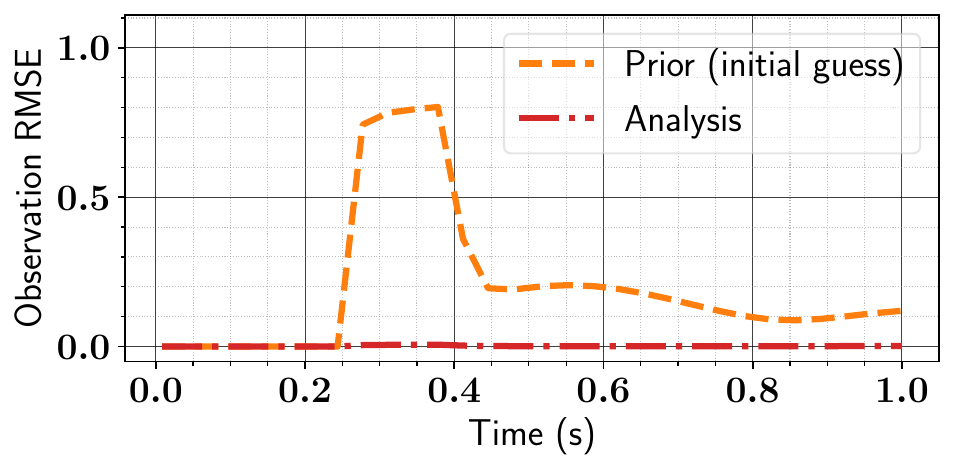}
            \includegraphics[trim={0 50pt 0 0},clip,width=0.75\linewidth]{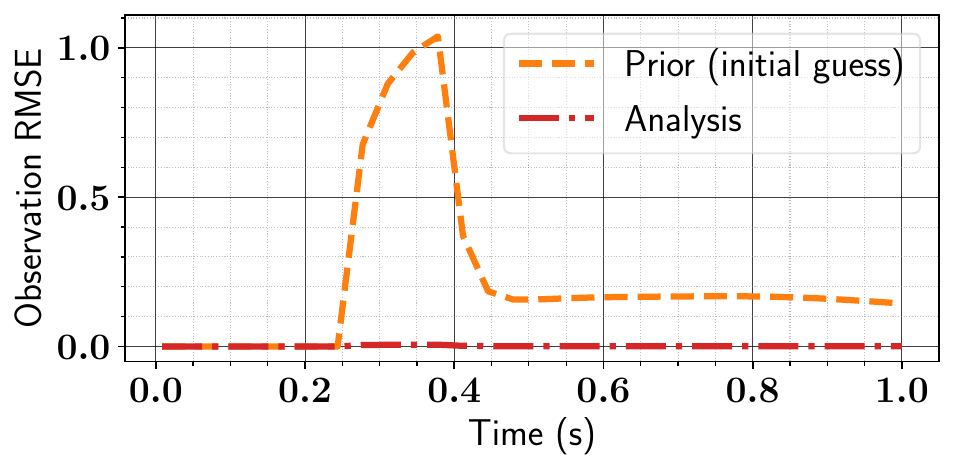}
            \caption{
                Inversion results obtained by applying Algorithm~\ref{alg:4DVar_Laplace_Approximation} with 
                fault impedance $0.01$. The fault occurs at buses $1,10,$ and $20$, respectively.
                For each choice of the faulty bus, the RMSE results are obtained by comparing the ground truth of the observations
                with model-based predicted observations at all observation time instances. 
                At each time instance, the RMSE is evaluated for all observations variables. 
            }
            \label{fig:nocorr-inversion-errors-fi-0.01-fb-1}
        \end{figure}

        Figure~\ref{fig:nocorr-inversion-errors-fi-0.01-fb-1} shows results of Algorithm~\ref{alg:4DVar_Laplace_Approximation} with the 
        fault impedance level set to $0.01$;  the fault occurs at buses $1$, $10$, and $20$, respectively.
        The small analysis RMSE results in the three cases show that the solution (analysis) obtained by applying \ref{alg:4DVar_Laplace_Approximation} 
        accurately estimates the unknown true parameter regardless of where the fault occurs in the network.

\begin{figure}[htbp!]
            \centering
            \includegraphics[trim={0 50pt 0 0},clip,width=0.75\linewidth]{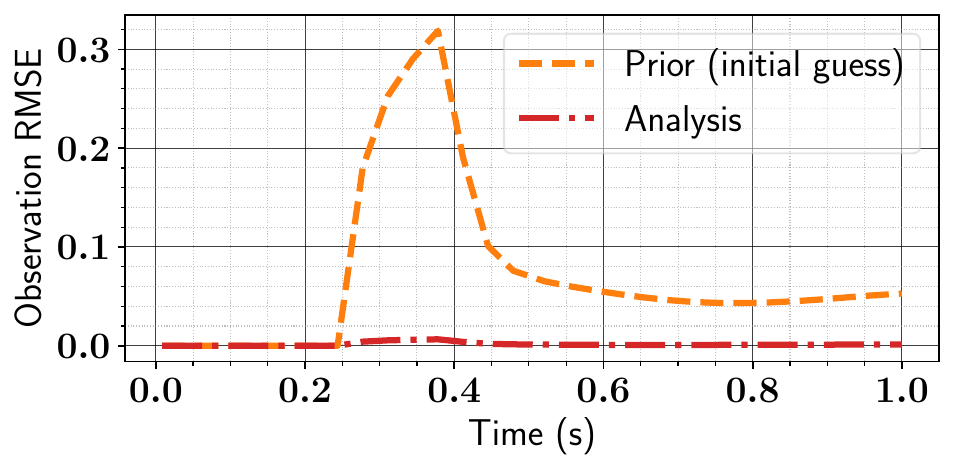}
            \includegraphics[trim={0 50pt 0 0},clip,width=0.75\linewidth]{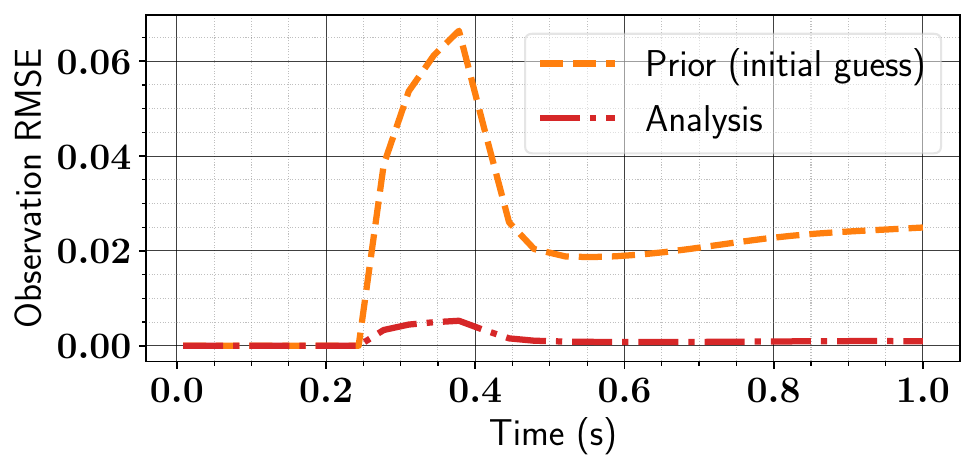}
            \includegraphics[trim={0 50pt 0 0},clip,width=0.75\linewidth]{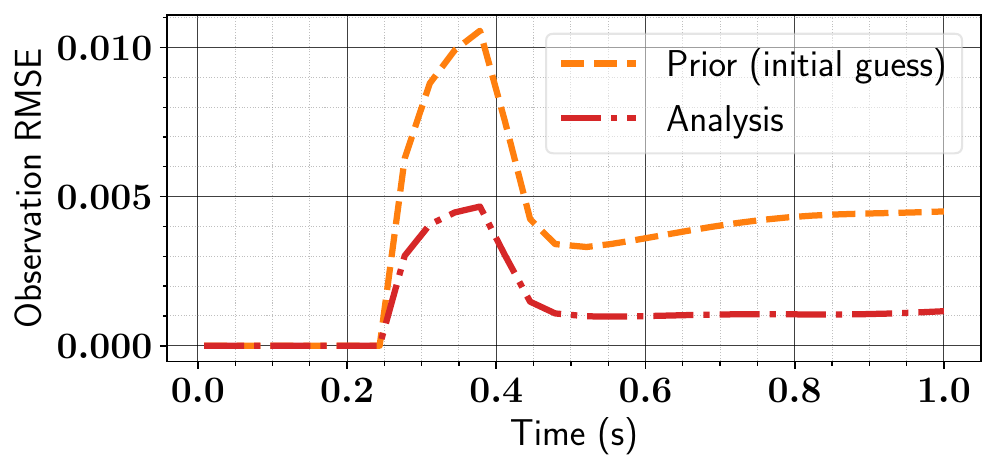}
            \caption{
                Similar to Figure~\ref{fig:nocorr-inversion-errors-fi-0.01-fb-1}.
                Here, the fault happens at bus $1$ and with the fault impedance level set to $0.03$, $0.1$, and $0.5$, respectively. 
            }
            \label{fig:nocorr-inversion-errors-fi-incr-fb-1}
        \end{figure}

\begin{figure}[htbp!]
            \centering
            \includegraphics[trim={0 50pt 0 0},clip,width=0.75\linewidth]{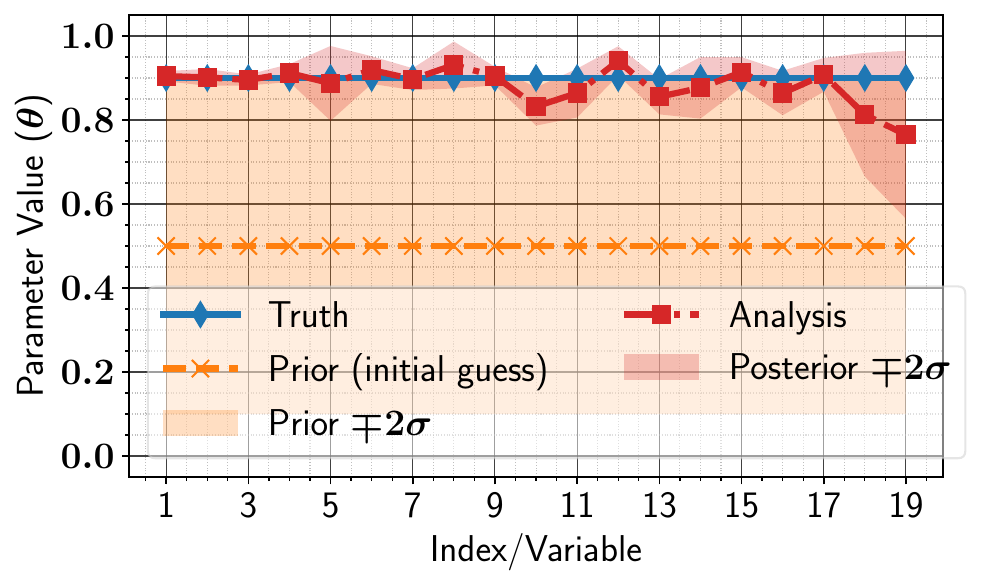}
            \includegraphics[trim={0 50pt 0 0},clip,width=0.75\linewidth]{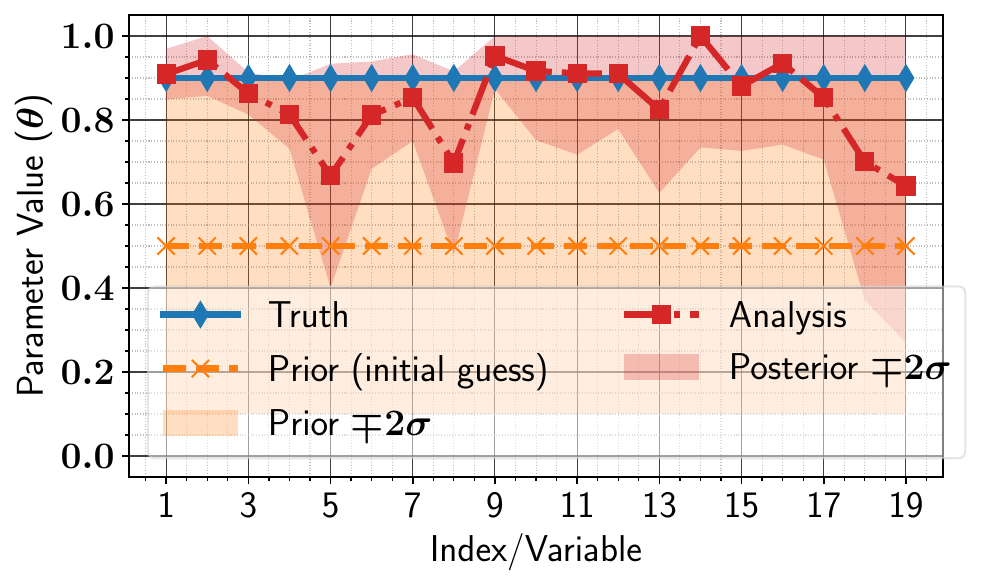}
            \includegraphics[width=0.75\linewidth]{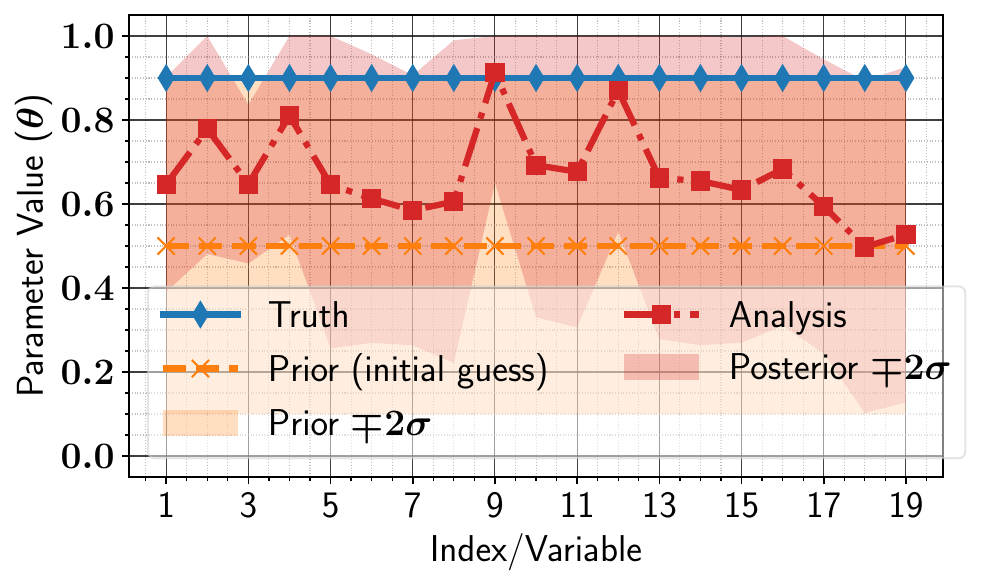}
            \caption{
                    Inversion results obtained by applying Algorithm~\ref{alg:4DVar_Laplace_Approximation} with fault at bus $1$ and with impedance level set to $0.03, 0.1, 0.5$, respectively.
                    Results here correspond to RMSE results shown in Figure~\ref{fig:nocorr-inversion-errors-fi-incr-fb-1}.
            }
            \label{fig:nocorr-inversion-fi-incr-fb-1}
        \end{figure}

        Figure~\ref{fig:nocorr-inversion-errors-fi-incr-fb-1} shows results of Algorithm~\ref{alg:4DVar_Laplace_Approximation} with the 
        fault synthesized at bus $1$, with fault impedance level set to $0.03$, $0.1$, and $0.5$, respectively.
        In all settings, the analysis RMSE is almost identical; however, with increasing impedance levels, the initial guess (prior mean) 
        yields a trajectory that is sufficiently accurate. This is expected because higher impedance results in smaller transients.
        Conversely, for lower impedance levels the need for inversion is evidently critical as the initial guess produces 
        observation trajectories with significantly higher RMSEs.
        The retrieved solution (inferred parameter) along with prediction uncertainty cones is shown in Figure~\ref{fig:nocorr-inversion-fi-incr-fb-1}, 
        for different contingencies. These results show that the inferred parameter becomes more erroneous as the impedance level increases; concurrently, the error in predicted observations is smaller because the transient is less acute. 

        Note that similar results were obtained by carrying out the same set of experiments with various choices of the failure location.
        Results shown in Figures~\ref{fig:nocorr-inversion-errors-fi-0.01-fb-1},~\ref{fig:nocorr-inversion-errors-fi-incr-fb-1}, 
        and~\ref{fig:nocorr-inversion-fi-incr-fb-1} indicate the capability of Algorithm~\ref{alg:4DVar_Laplace_Approximation} to accurately infer 
        the unknown parameter in any contingency and provide consistent UQ.

    \subsubsection{Predictivity analysis}
    \label{subsubsec:predectivity}
        We now analyze the predictive properties of the model calibrated by using Algorithm~\ref{alg:4DVar_Laplace_Approximation}. 
        We proceed by applying a fault of a certain impedance and bus. 
        We calibrate the model based on observations that result from the ensuing transient and then use the calibrated model to predict 
        transients that occur as a result of faults with different impedance levels and fault locations.

\begin{figure}[htbp!]
            \centering
            \includegraphics[width=0.85\linewidth]{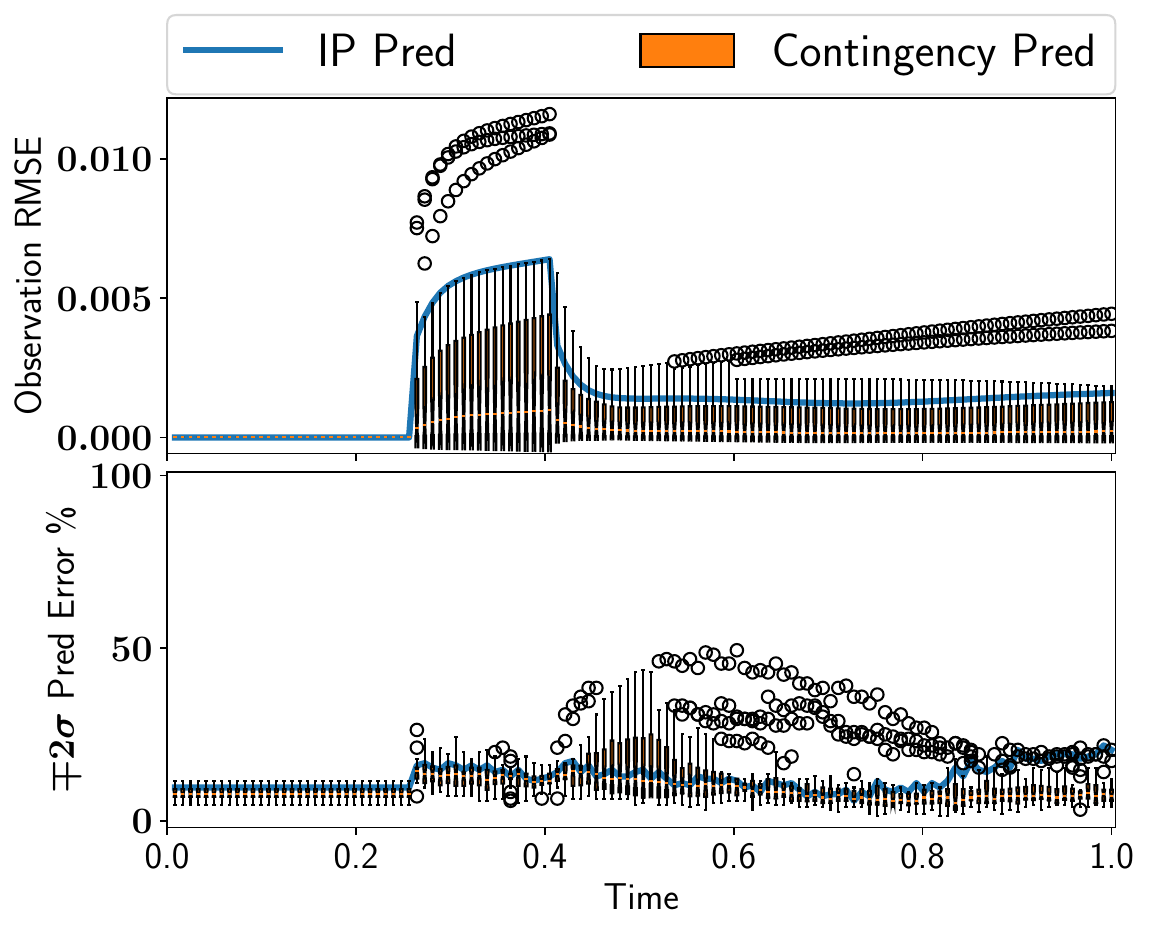}
            \caption{
                Prediction results  obtained by applying Algorithm~\ref{alg:4DVar_Laplace_Approximation} with fault impedance $0.01$/ The fault happens at bus $10$.
                The prediction is carried out for all  $20$ contingencies, and the results are aggregated at each time instance.
                The RMSE (top) is combined and plotted (box plot) at each time instance for all observation variables for all contingencies.
                The coverage error (bottom) is defined as ($100\%$ - coverage) where coverage  is the percentage of actual observations within $\mp 2\vec{\sigma}$ from the predicted observations, 
                where $\vec{\sigma}$ is the predicted/posterior standard deviation projected onto the observation space.
            }
            \label{fig:nocorr-prediction-errors-fi-0.01-fb-1}
        \end{figure}

        We noticed that (a) the results obtained by a specific impedance level are similar for all 
        choices of the faulty bus and (b) with increasing values of the impedance level, the inversion results become exceedingly erroneous. 
        Thus, for clarity, we set the faulty bus to $10$ and show results only for the impedance level
        $0.01$; see Figure~\ref{fig:nocorr-prediction-errors-fi-0.01-fb-1}. 

        Figure~\ref{fig:nocorr-inversion-fi-incr-fb-1} shows that the accuracy of the analysis (inferred parameter) is much better for lower 
        values of the impedance level. However, the posterior uncertainty (e.g., posterior variance) tends to be underestimated, resulting in overfitting, 
        which affects the predictive power in case of contingencies other than the one used for solving the inverse problem.
        This is explained by the results shown in Figure~\ref{fig:nocorr-prediction-errors-fi-0.01-fb-1}.

        %

\begin{figure}[!htbp]
            \centering
            \includegraphics[width=0.49\linewidth]{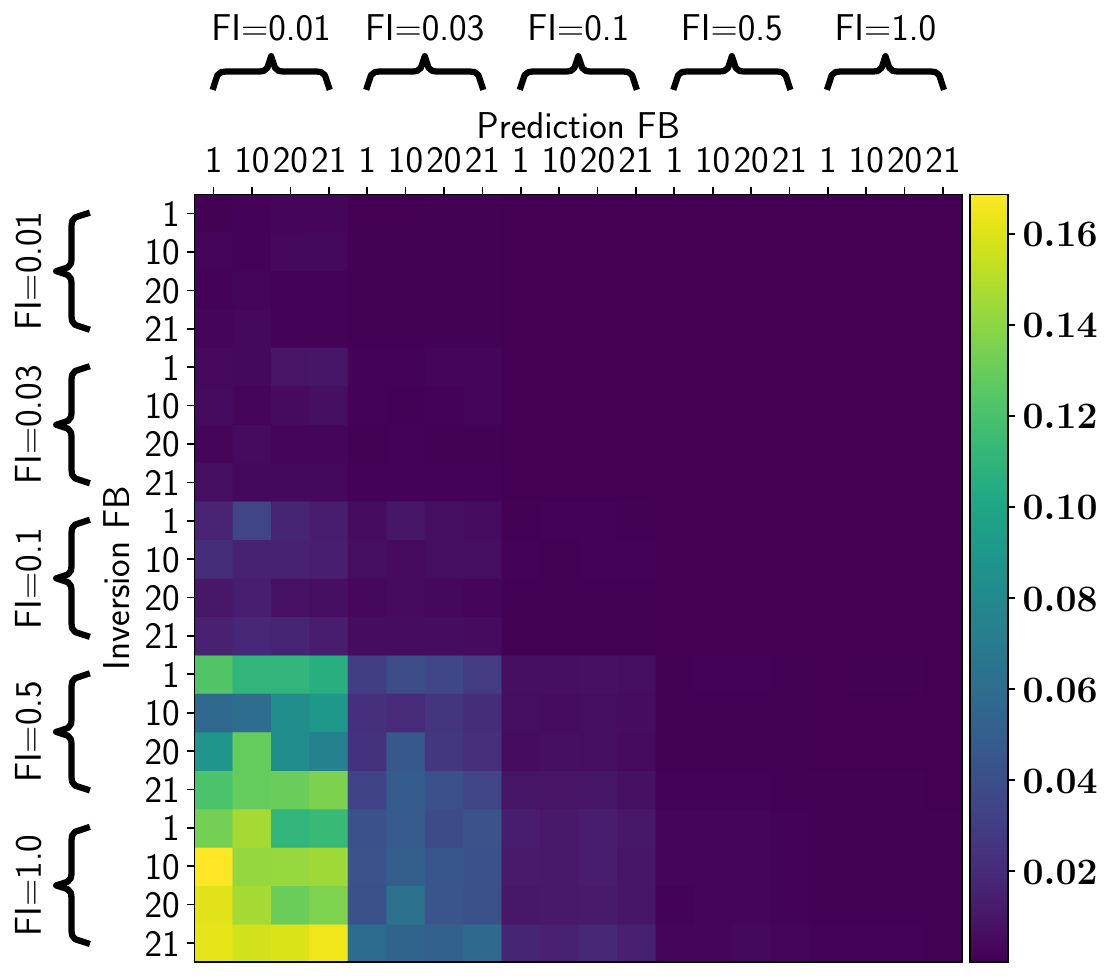}
            \includegraphics[width=0.49\linewidth]{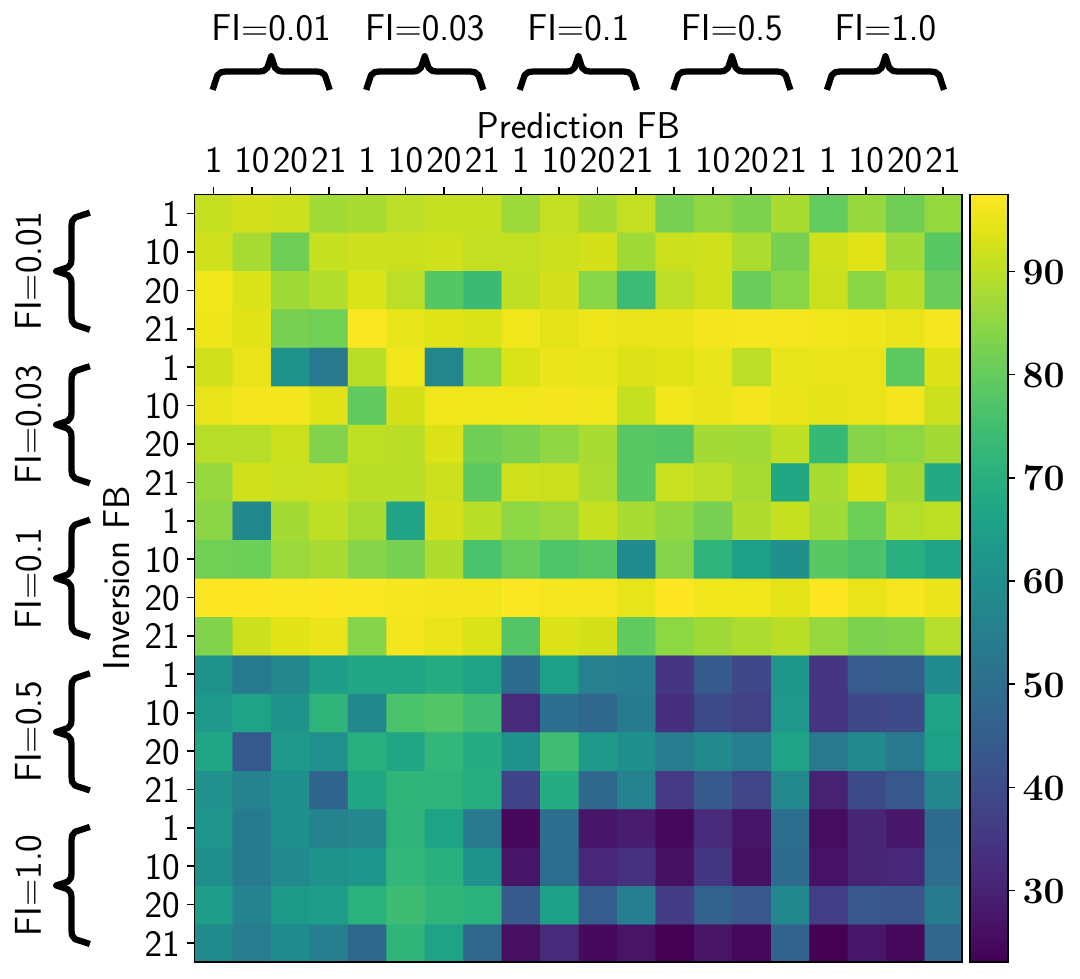}
            \caption{
                DA results  obtained by various choices of contingency, that is, the choice of the fault impedance level and the faulty bus.
                For each contingency, the solution obtained (MAP and posterior covariance) is used to make a prediction for each contingency. 
                Left: RMSE values (average is taken over all variables and all time instances in the assimilation window).
                Right: the coverage is defined as the percentage of actual observations within $\mp 2\vec{\sigma}$ from the predicted observations, 
                where $\vec{\sigma}$ is the predicted/posterior standard deviation projected onto the observation space.
                No covariance inflation.
            }
            \label{fig:nocorr-noinfl-prediction-rmse-coverage}
        \end{figure}

        Figure~\ref{fig:nocorr-prediction-errors-fi-0.01-fb-1} 
        shows that Algorithm~\ref{alg:4DVar_Laplace_Approximation}
        retrieves a good estimate of the true parameter (explained by the low levels of RMSE) 
        for both inversion and predicting other contingencies. 
        The prediction accuracy (for other contingencies than the one used for parameter identification) explained by the box plots fluctuates around 
        the inversion results in some relatively high fluctuations resulting from overfitting. 
        The effect of overfitting is further explained by Figure~\ref{fig:nocorr-noinfl-prediction-rmse-coverage}.

\begin{figure}[!htbp]
            \centering
            \includegraphics[width=0.49\linewidth]{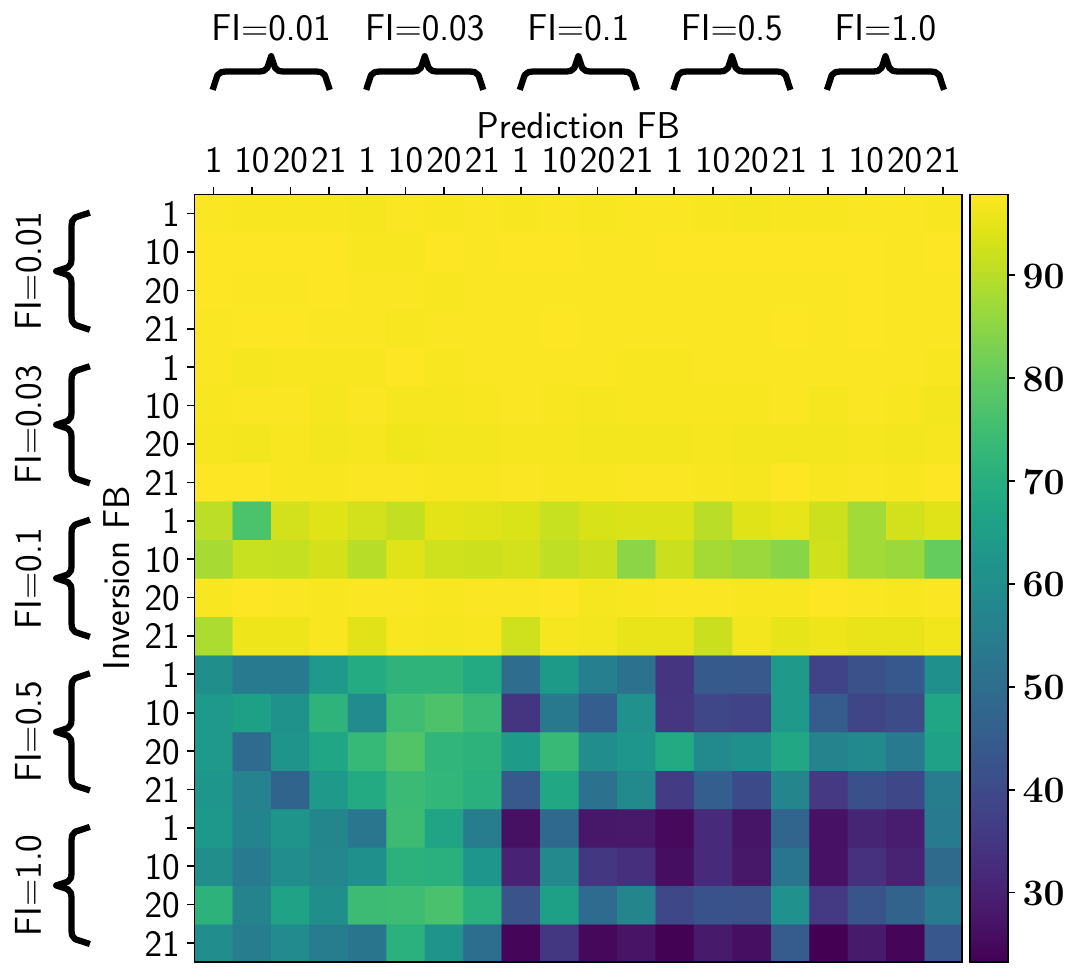}
            \includegraphics[width=0.49\linewidth]{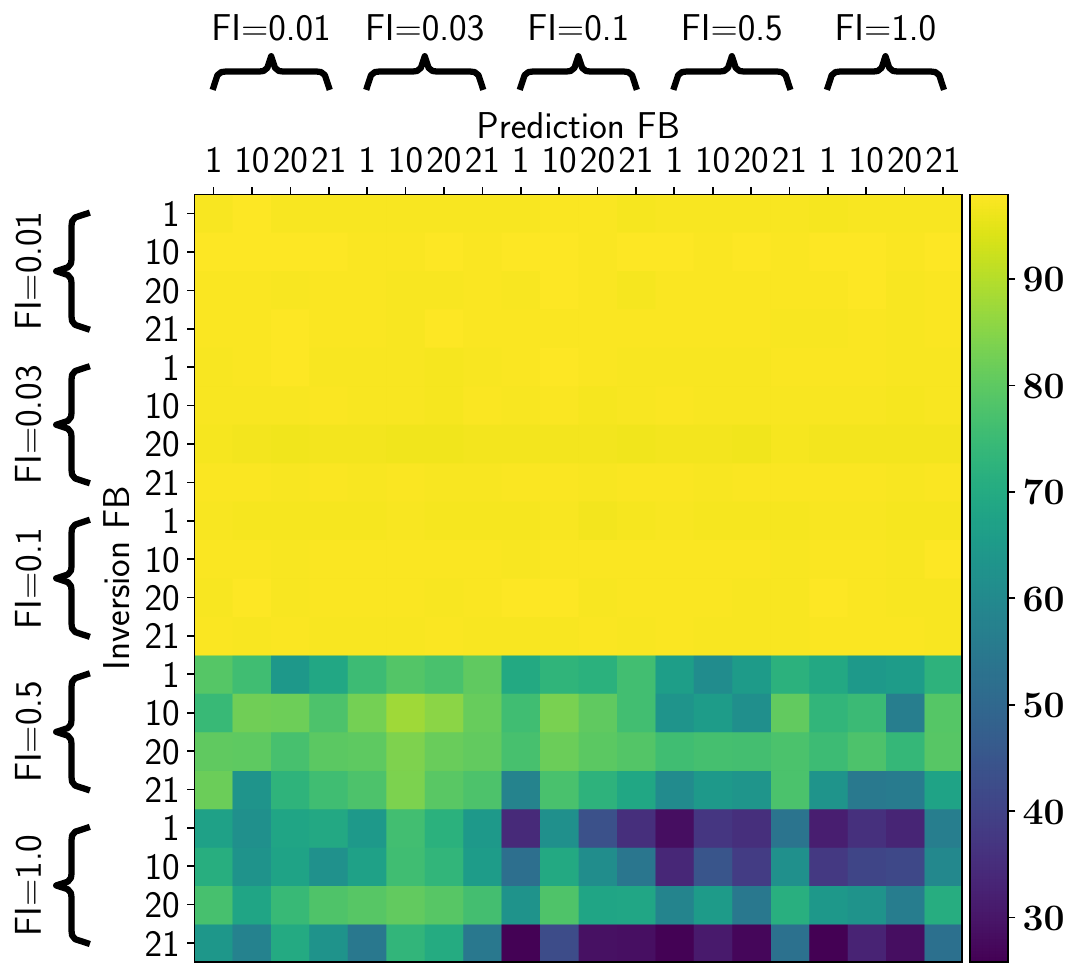}
            \caption{
                Similar to Figure~\ref{fig:nocorr-noinfl-prediction-rmse-coverage}. 
                Here we show coverage results obtained by applying various levels of covariance inflation to 
                the posterior covariance matrix.
                Results are shown for the inflation factor values $\lambda=0.001$ and $\lambda=0.005$, respectively. 
            }
            \label{fig:nocorr-infl-prediction-coverage}
        \end{figure}

        Overfitting can be ameliorated by applying covariance inflation~\cite{anderson2007adaptive,attia2018optimalinflation}, 
        a technique commonly used in the DA literature to cope with sampling errors resulting from employing 
        small sample sizes for covariance estimation.
        Here, covariance inflation is carried out by adding a constant scalar value $\lambda$ to the analysis (posterior) 
        variance; that is, the posterior covariance matrix $\postcov$ is replaced with $\postcov+\lambda \mat{I}$, where $\lambda > 0$ 
        is a positive inflation factor. Of course, inflation does not affect the solution accuracy because inflation is applied to 
        the posterior covariance, and it does not affect the MAP estimate. Inflation only enables wider coverage of observations 
        through prediction. 
        This is explained by results in Figure~\ref{fig:nocorr-infl-prediction-coverage}, which show that even small values of covariance
        inflation factor (e.g., $\lambda=0.001$) yield better coverage, that is, lower prediction errors.
        Note that for higher values of the impedance level, the coverage results are generally poor, which supports our previous assertion
        that inversion must be carried out with small values of the impedance levels in order to properly predict other contingencies.

\begin{figure}[!htbp]
            \centering
            \includegraphics[width=0.75\linewidth]{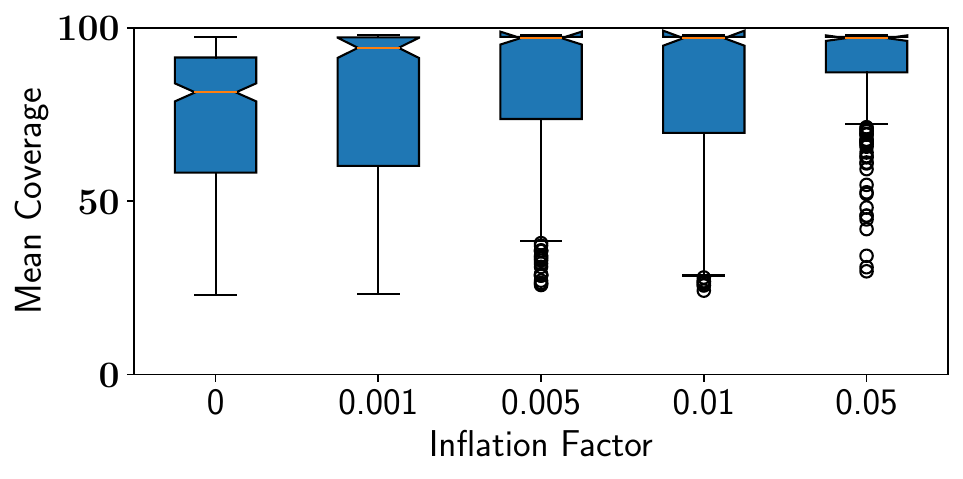}
            \caption{
                Coverage as a function of the inflation level.
                The coverage here is summarized by taking the average coverage over all contingencies 
                displayed in Figures~\ref{fig:nocorr-noinfl-prediction-rmse-coverage} and Figure~\ref{fig:nocorr-infl-prediction-coverage}. 
            }
            \label{fig:nocorr-coverage-boxplot}
        \end{figure}

        Figure~\ref{fig:nocorr-coverage-boxplot} further summarizes the results in Figure~\ref{fig:nocorr-infl-prediction-coverage} 
        with additional values of the inflation factor and shows the coverage results for increasing levels of the inflation factor $\lambda$.
        Note, however, that as the inflation factor increases, the uncertainty cones are expected to be wider. In this case, however, the
        modified covariance matrix (obtained by inflation) differs significantly from the posterior covariance obtained by Algorithm~\ref{alg:4DVar_Laplace_Approximation} with the posterior variances being overestimated,  thus degrading the amount of information
        provided by the UQ estimate.
        An analysis is required to obtain an optimal value of the inflation factor, which is left for future work.

    \subsubsection{TLM-based inversion and prediction}
    \label{subsubsec:tlm_results}
        We conclude this section by comparing the performance of Algorithm~\ref{alg:4DVar_Laplace_Approximation}
        with the case where the TLM is obtained beforehand and is used in the analysis (both inversion and prediction) 
        instead of the full nonlinear model. Doing so will allow us to evaluate the benefit of computing the posterior mode followed by the system linearization, as we propose here, compared with linearizing the system a priori at the reference point, as was carried out in \cite{nagi2021bayesian}.
        
\begin{figure}[!htbp]
            \centering
            \includegraphics[width=0.75\linewidth]{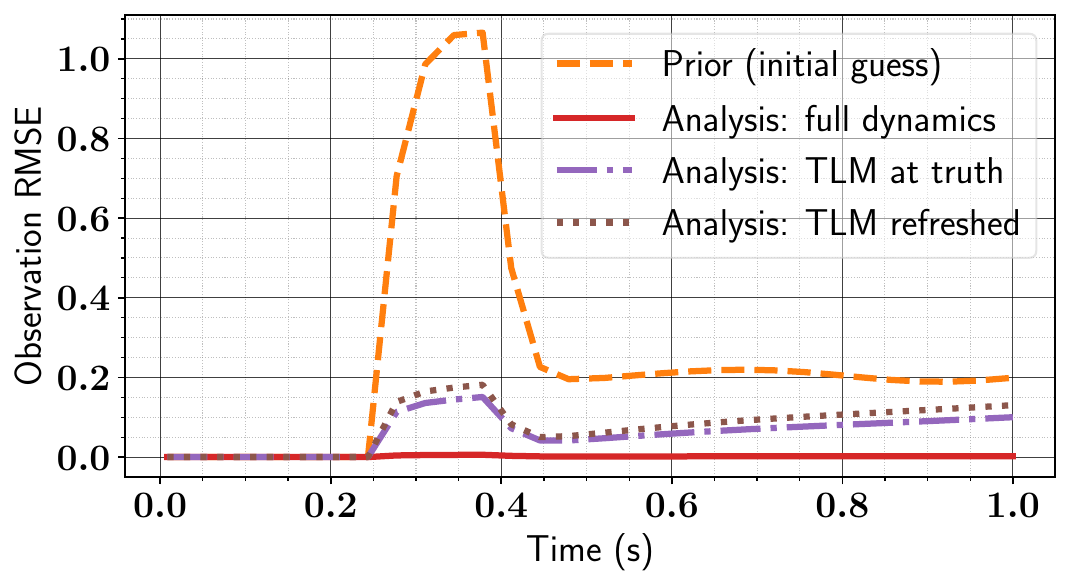}
            \includegraphics[width=0.75\linewidth]{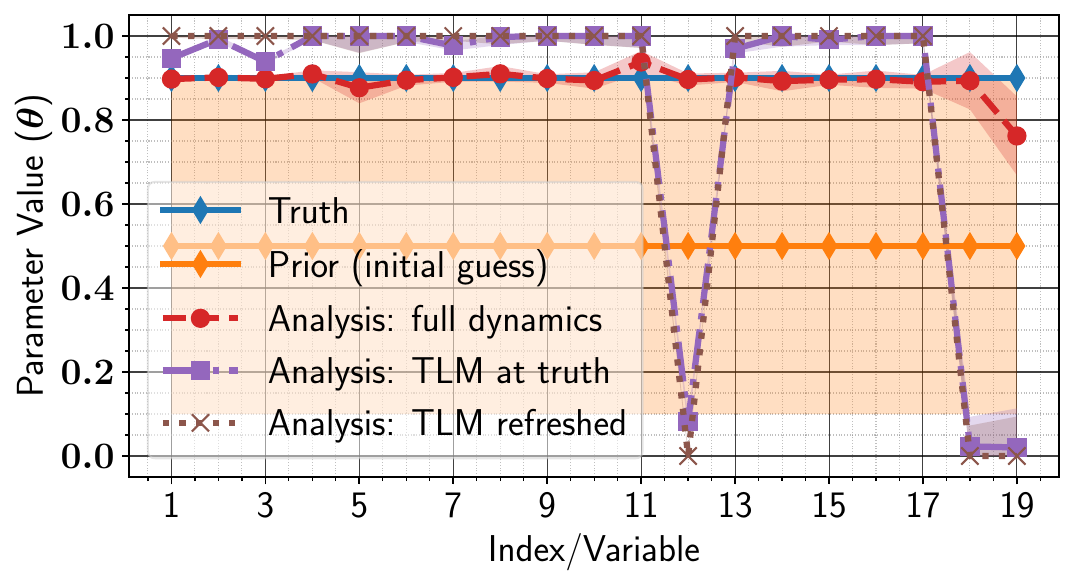}
            \caption{
                Inversion results obtained with Algorithm~\ref{alg:4DVar_Laplace_Approximation} 
                using full model dynamics compared with using the TLM for forward model simulations.
                The TLM  either is evaluated at the ground truth of the parameter or is reevaluated by linearizing model 
                dynamics around the parameter at which forward simulation is required.
                Results are obtained by setting the impedance level to $0.01$;  the fault happens at bus $1$.
                Top: RMSE results.
                Bottom: inversion parameter with uncertainty cones.
            }
            \label{fig:nocorr-TLM-RMSE}
        \end{figure}

        Here we show results for two cases and compare each with  Algorithm~\ref{alg:4DVar_Laplace_Approximation}; 
        see Figure~\ref{fig:nocorr-TLM-RMSE}.
        In the first case, we evaluate the TLM (forward sensitivities) using the ground truth, which is generally 
        unavailable in practice. This helps  create a benchmark for TLM involvement.
        In the second case,  the TLM is reevaluated/refreshed by linearizing model dynamics around the 
        parameter at which forward simulation is required, for example, when the gradient is evaluated in the optimization procedure. 
        
        The CPU times for the average cost of one objective, gradient, and linearization for both approaches are given in Table \ref{tab:CPU:time}.
        \begin{table}[!htbp]
            \begin{tabular}{|l|p{2cm}|c|c|}
            \hline
            Method & Objective \hspace{0.5cm} \quad (fwd. integration) & Gradient & Linearization\\
            \hline
            Full dynamics  & 1.03 (0.872) & 12.41 & - \\
            \hline
            TLM & 0.05 (0.001) & 0.10 &  12.60\\
            \hline
            TLM (refreshed) & 0.05 (0.001) & 12.70 &  12.60\\
            \hline
            \end{tabular}
            \caption{CPU time [seconds] estimation for the proposed method (full dynamics) and the linearized approach (TLM) averaged over 100 samples. The proposed nonlinear method is slower but more accurate. The linearized method is faster after computing and storing the linearized operators -- the linearization is done once if the operators can be stored in memory.} 
            \label{tab:CPU:time}
        \end{table}

        Results in Figure~\ref{fig:nocorr-TLM-RMSE} show that refreshing the TLM achieves results similar to the case 
        where the TLM is evaluated at the ground truth. In both cases, however, the solution is biased. 
        Algorithm~\ref{alg:4DVar_Laplace_Approximation}, on the other hand, outperforms the utilization of TLM with unbiased
        and significantly more accurate results.
        
       Comparing the total cost with the TLM method proposed in \cite{nagi2021bayesian} is more difficult. 
       The computational workflow of \cite{nagi2021bayesian}, required to determine the reference value at which to obtain the TLM approximation, is  complex and does not resemble the other computations we do here to produce a crisp statement of comparison. Our problem  has 19 parameters; and, on average, the nonlinear method we propose took 23 iterations to produce the MAP estimator and the quantities required to carry out the uncertainty quantification, which would be about 270 seconds overall with our Python framework.  If the point at which to compute the TLM were known, then our framework would require only about 12 seconds to carry out the computation, certainly much faster than our iterative nonlinear approach. On the other hand, the point is not known. Computing the proper approximation by  the method from \cite{nagi2021bayesian}, which is coded in MATLAB, was reported in that reference to take 180 seconds for 20 parameters for one iteration, with the total number of iterations unreported. While the overall computational expenditure is difficult to compare, we find that our algorithm is at least competitive with the method from \cite{nagi2021bayesian} in compute time and superior in accuracy as we reported above.


    \section{Conclusions}
    In this paper we have presented a novel centralized, Bayesian variational data assimilation technique to calibrate power system dynamic models. Our method is able to obtain posterior distribution estimates of parameter values and quantify the uncertainty of the inference. By experimenting with different fault locations and intensities, we show how our method produces uncertainty estimates that are consequent with the amount of information in the observed data. Furthermore, we show how our method is able to improve over previous state of the art based on the TLM approximation. 
For future work, it would be useful to profile the relative cost accuracy balance of other data assimilation methods, such as ensemble Kalman  or particle filters, since despite their known asymptotic disadvantages compared with variational methods 
they may of course do better on a problem at fixed size. Given our focus on comparing different workflows of variational data assimilation methods, we postpone such analysis for future research.

    \section*{Acknowledgments}
    This material was based upon work
    supported by the U.S. Department of Energy, Office of Science,
    Office of Advanced Scientific Computing Research (ASCR) under
    Contract DE-AC02-06CH11347. 
    
    
    \bibliographystyle{ieeetr}
    \bibliography{Bib/main,Bib/ahmed_references}
    
    \vspace{0.3cm}
    \begin{center}
        \scriptsize \framebox{
            \parbox{2.5in}{
                Government License (will be removed at publication):
                The submitted manuscript has been created by UChicago Argonne, LLC,
                Operator of Argonne National Laboratory (``Argonne").  Argonne, a
                U.S. Department of Energy Office of Science laboratory, is operated
                under Contract No. DE-AC02-06CH11357.  The U.S. Government retains for
                itself, and others acting on its behalf, a paid-up nonexclusive,
                irrevocable worldwide license in said article to reproduce, prepare
                derivative works, distribute copies to the public, and perform
                publicly and display publicly, by or on behalf of the Government. The Department of Energy will provide public access to these results of federally sponsored research in accordance with the DOE Public Access Plan. http://energy.gov/downloads/doe-public-access-plan.
            }
        }
        \normalsize
    \end{center}
    
\end{document}